\documentclass[12pt]{amsart}
\usepackage[margin=1in]{geometry}
\usepackage{ifxetex,ifluatex, tikz, pgfplots}
\usepackage{makecell}

\newif\ifxetexorluatex
\ifxetex
  \xetexorluatextrue
\else
  \ifluatex
    \xetexorluatextrue
  \else
    \xetexorluatexfalse
  \fi
\fi

\ifxetexorluatex
  \usepackage{fontspec}
\else
  \usepackage[T1]{fontenc}
  \usepackage[utf8]{inputenc}
  \usepackage{lmodern}
\fi
\usepackage{hyperref, amsmath, amssymb, mathtools}

\usepackage{amsthm}
\usepackage{amsfonts}
\usepackage{amsmath}
\usepackage{comment}
\usepackage{chngcntr}
\usepackage{ctable}
\newtheorem{thm}{Theorem}
\newtheorem*{thm*}{Theorem}

\newtheorem{corollary}[thm]{Corollary}

\theoremstyle{definition}

\newtheorem*{remark}{Remark}
\numberwithin{thm}{section}

\newcommand{\Z}{\mathbb{Z}}
\newcommand{\Q}{\mathbb{Q}}
\newcommand{\set}[1]{\left\{ #1 \right\}}
\renewcommand{\H}{\mathbb{H}}
\newcommand{\SL}{\text{SL}_2(\Z)}
\newcommand{\C}{\mathbb{C}}
\let\Re\relax
\DeclareMathOperator{\Re}{Re}
\renewcommand{\Im}{\text{Im}}

\title{On the zeros of a class of modular functions}
\author[N. Sweeting]{Naomi Sweeting}
\address{Department of Mathematics, University of Chicago, Chicago, IL 60637}
\email{nsweeting@uchicago.edu}
\author[K. Woo]{Katharine Woo}
\address{Department of Mathematics, Stanford University, Stanford, CA 94305}
\email{katywoo@stanford.edu}

\date{July 2018}

\begin{document}


\maketitle
\vspace*{-1cm}
\begin{abstract}
    We generalize a number of works on the zeros of certain level 1 modular forms to a class of weakly holomorphic modular functions whose $q$-expansions satisfy \[
    f_k(A; \tau) \coloneqq q^{-k}(1+a(1)q+a(2)q^2+...) + O(q),\] where $a(n)$ are numbers satisfying a certain analytic condition. We show that the zeros of such $f_k(\tau)$ in the fundamental domain of $\SL$ lie on $|\tau|=1$ and are transcendental. We recover as a special case earlier work of Witten on extremal ``partition'' functions $Z_k(\tau)$. These functions were originally conceived as possible generalizations of constructions in three-dimensional quantum gravity.
\end{abstract}

\section{Introduction and Statement of Results}
Weakly holomorphic modular functions on $\SL$ can be uniquely constructed to have arbitrary principal parts in their $q$-expansion. Given a formal power series $$A(q) \coloneqq 1 + a(1) q + a(2)q^2+ \ldots,$$ we define an associated modular function $f_k(A; \tau)$ as the unique weakly holomorphic modular function on $\SL$ satisfying 
\begin{equation}
f_k(A; \tau) \coloneqq q^{-k} A(q) + O(q).
\end{equation} (Throughout, $q \coloneqq e^{2\pi i \tau}$.) Using Faber polynomials $F_n(j)$ (see (\ref{faber-poly})), which are monic degree $n$ polynomials in the $j$ invariant (see (\ref{j})), we can write $f_k(A; \tau)$ as a degree $k$ polynomial in $j(\tau)$:
\begin{equation}\label{poly-expand}
f_k(A; \tau) = P_k(A; j(\tau)) \coloneqq F_k(j(\tau)) + \sum_{n=0}^{k-1} a(k-n)F_n(j(\tau)).
\end{equation}

We investigate the zeros of $f_k(A; \tau)$ by extending the results of Rankin and Swinnerton-Dyer \cite{RSD} and Asai, Kaneko, and Ninomiya \cite{AKN}, who studied the zeros of Eisenstein series and Faber polynomials, respectively. To state our result, we use a bound from \cite{AKN}. For $\tau \in C\coloneqq  \set{\tau \in \H \, \colon |\tau| = 1, 0 \leq \Re(\tau) \leq 1/2}$, we have:
\begin{align}\label{faber-zeros}
    |F_n(j(\tau)) e^{-2\pi n \Im(\tau)} - 2\cos(2\pi n\Re(\tau))| &< M \approx 1.1176\ldots 
\end{align} 
This bound is useful because $j(C) = [0, 1728] \subset \mathbb{R}$, so the intermediate value theorem allows us to identify at least one zero in each of the $k$ segments \[C_{n,k} \coloneqq C \cap \set{\frac{n-1}{2k} < \Re(x) < \frac{n}{2k}}, \;  n = 1, \ldots, k.\] Using this notation, we now state our first theorem. 
\begin{thm}
Let $A(q)$, $f_k(A; q)$, and $a(n)$ be as above. Suppose that
\[
S\coloneqq \sum_{n = 1}^\infty |a(n)| e^{-\pi n\sqrt{3}}  < \frac{2-M}{2+M}.
\]  Then, for all integers $k \geq 1$, $f_k(A; \tau)$ has exactly $k$ zeros in the fundamental domain for $\SL$, and they all lie on $C$. Moreover, each segment $C_{n, k}$ of the arc contains exactly one zero of $f_k(\tau).$
\end{thm}

\begin{remark}
This result gives equidistribution of the zeros of $f_k (A; \tau)$ with respect to $\Re(\tau)$. 
\end{remark}

Because $f_k(A; \tau) = P_k(j(A; \tau))$ and $j$ is a bijection between the fundamental domain and $\C$, our theorem yields the following corollary.
\begin{corollary}
Under the hypotheses of Theorem 1.1, $P_k(A; y)$ has $k$ distinct real zeros in the interval $(0, 1728).$
\end{corollary}

In 2007, Witten \cite{Witten} defined  functions that conjecturally encode the existence of extremal conformal field theory at any central charge $c = 24k$ with $k\geq 1$. He suggested that for $k\geq 1$, the generating functions for the quantum states of three dimensional gravity are given by
\[
Z_k(\tau) \coloneqq q^{-k} \prod_{n\geq 2} (1-q^n)^{-1} + O(q).
\] 
Theorem 1.1 and Corollary 1.2 apply to all of these $Z_k$, as illustrated in Figure 1. Thus we have recovered the following result of Witten \cite{WittenM}. \begin{corollary}
For $k\geq 1$, $Z_k(\tau)$ has exactly $k$ zeros in the fundamental domain for $\SL$ and they all lie on $C$. In addition, each $C_{n,k}$ contains exactly one zero of $Z_k(\tau)$.
\end{corollary}

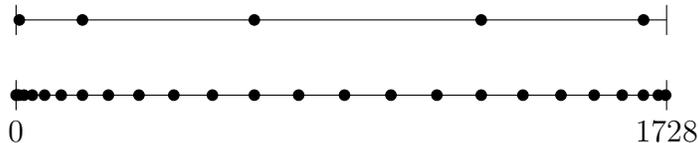
\begin{figure}[h!]
\begin{center}
\begin{tikzpicture}[only marks, y=1cm, x = 0.005cm]
    \draw plot[mark=*] file {Z5roots.data};
    \draw plot[mark=*] file {Z25roots.data};
    \draw (0,0) -- (1728,0);
    \draw (0, 1) -- (1728, 1);
    \draw (0, 0.2) -- (0, -0.2) node[anchor = north] {$0$};
    \draw (1728, 0.2) -- (1728, -0.2) node[anchor = north] {$1728$};
    \draw (0, 0.8) -- (0, 1.2);
    \draw (1728, 0.8) -- (1728, 1.2);
\end{tikzpicture}

\caption{The roots of the polynomials for $Z_5$ and $Z_{25}$.}
\end{center}
\end{figure}

Kohnen \cite{Kohnen} proved that all of the zeros except $i, e^{2\pi i/3}$ of the Eisenstein series and $F_k(j(\tau))$ are transcendental. We show a similar result for $f_k(A; \tau)$.

\begin{thm}
Let $f_k(A; \tau)$ be as above, and suppose that $a(n) \in \Q$ satisfy the hypotheses of Theorem 1.1. If $f_k(A; \tau_0) = 0$, then $\tau_0$ is transcendental. 
\end{thm}

\section{Preliminaries}\label{prelims}
First we briefly discuss modular functions; a general survey of these and other facts from the theory of modular forms can be found in \cite{Ono}. A modular function is a meromorphic modular form of weight zero on $\SL$. There are no nonconstant holomorphic modular functions, but we call a modular function weakly holomorphic if it has no poles on $\H$. 

The first example is the $j$-function, a weakly holomorphic modular function defined by \begin{equation}\label{j}j(\tau) \coloneqq  \frac{\left(1 + 240\sum_n \sigma_3(n)q^n\right)^3}{q \prod_n (1 - q^n)^{24}} = q^{-1} + 744 + 196884q + \ldots. \end{equation} If we define the fundamental domain for the action of $SL_2(\Z)$ on $\H$ by \begin{equation}\mathcal{F} \coloneqq \set{\tau \in \H\, \colon |\tau| \geq 1, -1/2 < \Re(\tau) \leq 1/2, -1/2 < \Re(\tau) < 0 \implies |\tau| > 1},\end{equation} then $j:\mathcal{F} \to \C$ is bijective.
All modular functions on $\SL$ are given by rational functions in $j$; the weakly holomorphic ones are the polynomials in $j$.

Any weakly holomorphic modular function is fully determined by the principal part of its $q$-expansion. We therefore define the Faber polynomial  $F_k(y)$, for an integer $k \geq 0$, to be the unique polynomial satisfying \begin{equation}\label{faber cusp} F_k(j(\tau))\coloneqq q^{-k} + O(q).\end{equation}  These polynomials are given \cite{Ono-Rolen} by the denominator formula for the monster Lie algebra: \begin{align}\label{faber-poly}\sum_{k = 0}^\infty F_k(y) q^k &= \frac{1 - 24\sum_n \sigma_{13}(n)q^n}{q\prod_n (1 - q^n)^{24}}\cdot \frac{1}{j(\tau) - 744 - y} \\ \nonumber &= 1 + (y-744)q + (y^2-1488y + 159768)q^2+\ldots.\end{align} The denominator formula plays a key role in the theory of monstrous moonshine. For our purposes, it suffices to note that $F_k$ are monic degree $k$ polynomials with integer coefficients and that they satisfy (\ref{faber cusp}). 

\section{Proofs}
\subsection{Proof of Theorem 1.1}
Without loss of generality, we assume $a(n) \in \mathbb{R}$; if not, split $A$ into real and imaginary parts. It suffices to show that \[\big|f_k(A; \tau)e^{-2\pi k \Im (\tau)} - 2\cos(2\pi k \Re(\tau))\big| < 2 \] for all $\tau \in C.$  We define an error term for $\tau \in C$ by \[R_k(\tau) \coloneqq \big|f_k (A; \tau) e^{-2\pi k \Im \tau} - 2\cos(2\pi k \Re(\tau))\big|\] and  use (\ref{poly-expand}) and (\ref{faber-zeros}) to calculate an upper bound as follows:
\begin{align*}
R_k(\tau) &< \big|F_k(j)e^{-2\pi k \Im \tau} - 2\cos(2\pi k \Re(\tau))\big| 
 + \sum_{n = 0}^{k-1} |a(k-n)F_n(j)|e^{-2\pi k \Im \tau} \\
&< M +  
\sum_{n = 0}^{k-1} |a(k-n) F_n(j)| e^{-2\pi n \Im \tau} e^{-2\pi (k-n) \Im \tau} \\
&< M + (2 + M) \sum_{n = 1}^{k} |a(n)| e^{-\pi n \sqrt{3}} \\
&\leq M + (2 + M) S.
\end{align*} For $a(n)$ satisfying our hypotheses, we have $R_k(\tau) < 2$ for all $\tau \in C$. Thus there are at least $k$ zeros of $f_k(A; \tau)$ on $C$, one on each segment $C_{n, k}$ of the arc. Since $f_k(A; \tau) = P_k(A; j(\tau))$,  $P_k$ is a real polynomial of degree $k$, and $j$ is a bijection between $\mathcal{F}$ and $\C$, we conclude that there are exactly $k$ zeros of $f_k(A; \tau)$ in $\mathcal{F}$ and they are the ones lying on $C$.

\subsection{Proof of Corollary 1.3}
Ono and Rolen \cite{Ono-Rolen} gave  formulas for $Z_k(\tau)$ in terms of the partition function $p(n)$ and Faber polynomials $F_n(y)$. In particular, $a(n) = p(n) - p(n-1)$ for all $n \geq 1$. 
Kane \cite[Remark, p. 15]{Kane}  proved an explicit version of Ramanujan's asymptotic for $p(n)$: 
    \begin{equation}\label{partitionbound}
p(n) < \frac{5.5}{n}\exp\set{\pi\sqrt{2n/3}}.
\end{equation} The trivial estimates $a(n) < p(n)$ and the observation $a(1) = 0$ then suffice to show that the conditions of the theorem are satisfied. 
\qed

\subsection{Proof of Theorem 1.4}
Since $F_n(y)$ have integer coefficients and $a(n)$ are all rational, $P_k(A; y) \in \Q[y]$. Suppose  $f_k(A; \tau_0) = 0$ for some $\tau_0 \in \mathcal{F}$. Then $j(\tau_0)$ is a root of $P_k$ and hence algebraic. By the celebrated theorem of Schneider \cite{Schneider}, $j(\tau_0)$ is either transcendental or imaginary quadratic. If $j(\tau_0)$ is imaginary quadratic, then $\tau_0$ corresponds to a class of quadratic forms with discriminant $D<0$. Using the theory of complex multiplication and the fact that the roots of $P_k$ are conjugate over $\Q(\sqrt{D})$, we have that $P_k(A; j(\tau_1)) = 0$ for 
\[
\tau_1 = \begin{cases}\frac{i\sqrt{|D|}}{2}, & D\equiv 0 \pmod{4}, \\ \frac{-1+i\sqrt{|D|}}{2}, & D\equiv 1 \pmod{4}.\end{cases}
\]
However, we also have $|\tau_1| = 1$. Thus $D$ is either $-3$ or $-4$. We observe that $i, e^{2\pi i/3}$ are not roots of $f_k(A; \tau)$ since $|2\cos(2\pi k \Re(\tau))| = 2$ at these points. Therefore $\tau_0$ is transcendental. \qed 

\bibliographystyle{abbrv}
\bibliography{biblio}

\end{document}